\documentclass[12pt]{amsart}
\usepackage{amssymb,amsmath,amsfonts,latexsym,setspace}
\usepackage{bm}
\newcommand{\newsection}[1]{\setcounter{equation}{0} \section{#1}}
\newcommand{\bea}{\begin{eqnarray}}
\newcommand{\eea}{\end{eqnarray}}

\newcommand{\clb}{\mathcal{B}}

\newcommand{\cld}{\mathcal{D}}
\newcommand{\cle}{\mathcal{E}}

\newcommand{\clh}{\mathcal{H}}
\newcommand{\clk}{\mathcal{K}}
\newcommand{\cll}{\mathcal{L}}

\newcommand{\clq}{\mathcal{Q}}

\newcommand{\cls}{\mathcal{S}}

\newcommand{\z}{\bm{z}}
\newcommand{\w}{\bm{w}}

\newcommand{\raro}{\rightarrow}

\def \qed {\hfill \vrule height6pt width 6pt depth 0pt}
\def\textmatrix#1&#2\\#3&#4\\{\bigl({#1 \atop #3}\ {#2 \atop #4}\bigr)}
\def\dispmatrix#1&#2\\#3&#4\\{\left({#1 \atop #3}\ {#2 \atop #4}\right)}
\newcommand{\be}{\begin{equation}}
\newcommand{\ee}{\end{equation}}
\newcommand{\ben}{\begin{eqnarray*}}
\newcommand{\een}{\end{eqnarray*}}

\newcommand{\NI}{\noindent}

\newcommand{\bi}{\begin{itemize}}
\newcommand{\ei}{\end{itemize}}

\newtheorem{Theorem}{\sc Theorem}[section]
\newtheorem{Lemma}[Theorem]{\sc Lemma}
\newtheorem{Proposition}[Theorem]{\sc Proposition}
\newtheorem{Corollary}[Theorem]{\sc Corollary}
\newtheorem{Definition}[Theorem]{\sc Definition}
\newtheorem{Example}[Theorem]{\sc Example}
\newtheorem{Remark}[Theorem]{\sc Remark}
\newtheorem{Note}[Theorem]{\sc Note}
\newtheorem{Question}{\sc Question}
\newtheorem{ass}[Theorem]{\sc Assumption}
\newcommand{\bt}{\begin{Theorem}}
\def\beginlem{\begin{Lemma}}
\def\beginprop{\begin{Proposition}}
\def\begincor{\begin{Corollary}}
\def\begindef{\begin{Definition}}
\def\beginexamp{\begin{Example}}
\def\beginrem{\begin{Remark}}
\def\beginq{\begin{Question}}
\def\beginass{\begin{ass}}
\def\beginnote{\begin{Note}}
\newcommand{\et}{\end{Theorem}}
\def\endlem{\end{Lemma}}
\def\endprop{\end{Proposition}}
\def\endcor{\end{Corollary}}
\def\enddef{\end{Definition}}
\def\endexamp{\end{Example}}
\def\endrem{\end{Remark}}
\def\endq{\end{Question}}
\def\endass{\end{ass}}
\def\endnote{\end{Note}}

\begin{document}

\title{Analytic Model of Doubly Commuting Contractions}

\author[Bhattacharyya]{T. Bhattacharyya}
\address{Department of Mathematics, Indian Institute of Science, Bangalore 560012, India}
\email{tirtha@member.ams.org}

\author[Narayanan]{E. K. Narayanan}
\address{Department of Mathematics, Indian Institute of Science, Bangalore 560012,
India} \email{naru@math.iisc.ernet.in}

\author[Sarkar]{Jaydeb Sarkar}
\address{Indian Statistical Institute, Statistics and Mathematics Unit, 8th Mile, Mysore Road, Bangalore, 560059, India}
\email{jay@isibang.ac.in, jaydeb@gmail.com}
%\date{\today}


%\today


\thanks{TB and EKN are supported by Department of Science and
Technology, India through the project numbered SR/S4/MS:766/12 and
University Grants Commission, India via DSA-SAP. JS is partially supported by
NBHM, India through the research grant NBHM/R.P.64/2014.}

\maketitle

\begin{abstract}
An $n$-tuple ($n \geq 2$), $T = (T_1, \ldots, T_n)$, of commuting
bounded linear operators on a Hilbert space $\clh$ is doubly
commuting if $T_i T_j^* = T_j^* T_i$ for all $1 \leq i < j \leq
n$. If in addition, each $T_i \in C_{\cdot 0}$, then we say that
$T$ is a doubly commuting pure tuple. In this paper we prove that
a doubly commuting pure tuple $T$ can be dilated to a tuple of
shift operators on some suitable vector-valued Hardy space
$H^2_{\cld_{T^*}}(\mathbb{D}^n)$. As a consequence of the dilation
theorem, we prove that there exists a closed subspace $\cls_T$ of
the form
\[\cls_{T} := \sum_{i=1}^n \Phi_{T_i}
H^2_{\cle_{T_i}}(\mathbb{D}^n),\] where
%Based on this dilation result, we prove that for a doubly
%commuting pure tuple $T$ on $\clh$, there exists Hilbert spaces
$\{\cle_{T_i}\}_{i=1}^n$ are Hilbert spaces, $\Phi_{T_i} \in
H^\infty_{\clb(\cle_{T_i}, \cld_{T^*})}(\mathbb{D}^n)$ such that
each $\Phi_{T_i}$ ($1 \leq i \leq n$) is either a one variable
inner function in $z_i$, or the zero function. Moreover, $\clh
\cong \cls_T^\perp$ and
\[(T_1, \ldots, T_n) \cong P_{\cls_T^\perp} (M_{z_1}, \ldots,
M_{z_n})|_{\cls_T^\perp}.\]

\end{abstract}

\section{Introduction}

Consider a complex separable Hilbert space $\cle$ and a closed
subspace $\cls$ of $H^2_{\cle}(\mathbb D)$ that is invariant under
the operator $M_z$ on $H^2_{\cle}(\mathbb D)$, i.e.,
$$ M_z \cls \subseteq \cls.$$
Clearly, $T = P_{\cls^\perp} M_z |_{\cls^\perp}$ is a contraction.
But, moreover, $T^{*m}$ converges to $0$ strongly as $m
\rightarrow \infty$. This is the so called $C_{\cdot 0}$ property that
$T$ inherits from $M_z$.

In their pioneering work in the late 1960's, Sz.-Nagy and Foias
showed that for a contraction to qualify as $C_{\cdot 0}$, it must
be of the above form. See \cite{SNF}. More precisely, if $T$ is a
$C_{\cdot 0}$ contraction on a Hilbert space $\clh$, then there is
an $\cle$ as above and a subspace $\cls_T$ of $H^2_{\cle}(\mathbb
D)$ such that $\cls_T$ is invariant under $M_z$ and $T$ is
unitarily equivalent to $P_{{\cls_T}^\perp} M_z
|_{{\cls_T}^\perp}$. Here $\cle$ is explicit. Indeed, if we denote
by $ D_{T^*}$ the defect operator $( I - TT^* )^{1/2}$, then
$\cle$ is nothing but $\cld_{T^*}$, the closure of the range of
$D_{T^*}$. This result was just one part of the revelation. The
technique through which it was achieved was equally revealing.
They produced $\cls_T$ as the range of the multiplier
$M_{\theta_T}$ where $\theta_T$ is the characteristic function of
$T$. Thus, they gave a Beurling-Lax-Halmos form of $\cls_T$.

Recall that the characteristic function of a contraction $T
\in \clb(\clh)$ is defined by
\begin{equation}\label{chfn}\theta_T(z) = -T + D_{T^*} (I -
zT^*)^{-1} z D_T. \quad \quad (z \in \mathbb{D})\end{equation} We
refer to \cite{SNF} for more properties of this function.

Such an elegant characterization of all $C_{\cdot 0}$ contractions
obviously led to a search for such a phenomenon in the polydisk
and the Euclidean unit ball. The challenges in a several variables
situation are manifold. One first had to identify  the space that
would play the role of the Hardy space. For the ball, it became
clear only in the 1990's with works of Drury \cite{Dru}, Pott
\cite{Pott}, Popescu \cite{Pop-Poisson} and Arveson \cite{sub3}
that the natural space for this purpose on the Euclidean unit ball
is the one with reproducing kernel $\frac{1}{1 - \langle z , w
\rangle}$. It was shown in \cite{BES1} that the above mentioned
result of Sz.-Nagy and Foias can be generalized to the Euclidean
unit ball.

The case of the polydisk is more interesting. There is no
generalization of the Sz.-Nagy Foias result mentioned above to
this situation. There are invariant subspace results though due to
Ahern and Clark \cite{AC}, Mandrekar \cite{Man}, Rudin \cite{R}
and Izuchi, Nakazi and Seto \cite{INS}. As far as the model theory
results are concerned, there is a general framework due to
Ambrozie, Englis and Muller \cite{AEM}. They do have a
generalization of the $C_{\cdot 0}$ condition which although pretty
natural when stated in an abstract setting, is quite intractable
after specializing to the polydisk.

This brings us to what we are doing in this note. We consider a
commuting tuple of contractions $T = (T_1, T_2, \ldots , T_n)$
such that $T_i^* T_j = T_jT_i^*$ for $i \neq j$ (double
commutativity) and $T_i^{*m} \rightarrow 0$ strongly for each $i$.
Under these assumptions, we give an interesting generalization of
the Sz.-Nagy Foias result involving characteristic functions of
the individual contractions. En route, we produce a new proof of
the model.

The paper is organized as follows. In section 2, we review and
collect some of the preliminary concepts that will be useful. In
section 3, we obtain a dilation result for pure doubly commuting
tuple of contractions.  In section 4, we obtain a functional model
for the class of pure doubly commuting tuples of contractions. In
the final section, section 5, we establish a relationship between
the class of pure doubly commuting tuples of contractions and one
variable inner functions defined on the unit polydisc.

\newsection{Preliminaries}

Before we introduce a tuple of doubly commuting contractions, let
us briefly review the case of a single contraction $T \in
\clb(\clh)$ which is  $C_{\cdot 0}.$ Consider the vector valued
Hardy space $H^2_{\cld_{T^*}}(\mathbb D).$ The contraction $T$ is
then realized as $ P_{\clq_T} M_z|_{\clq_T},$ where $\clq_T$ is
the orthogonal complement of $M_{\theta_T} H^2_{\cld_{T}}(\mathbb
D).$ A key ingredient in this theory is the map $L_T : \clh \raro
H^2_{\cld_{T^*}}(\mathbb{D})$ defined by
\begin{equation}\label{e-1} L_T h := D_{T^*} (I - z T^*)^{-1} h =
\sum_{n=0}^\infty z^n D_{T^*} T^{*n} h. \quad \quad (h \in
\clh).\end{equation} Then $L_T$ is an isometry and
\begin{equation}\label{e-2}L_T T^* = M_z^*L_T.\end{equation}
Moreover,\begin{equation*}\label{L-1}L_T^* (\mathbb{S}_w \otimes
\eta) = (I - \bar{w} T)^{-1} D_{T^*} \eta,\quad \quad (w \in
\mathbb{D}, \eta \in \cld_{T^*})\end{equation*} and
\begin{equation}\label{e-3}\mathbb{S}(\lambda, w)(I - \theta_T(\lambda)
\theta_T(w)^*) = D_{T^*} (I - \lambda T^*)^{-1} (I - \bar{w}
T)^{-1} D_{T^*}, \quad \quad (\lambda, w \in
\mathbb{D})\end{equation} where $\mathbb{S}$ is the Szego kernel
on the unit disk defined by $\mathbb{S}(z,w) = (1 - z
\bar{w})^{-1}$ for all $z, w \in \mathbb{D}$. \NI The above two
equalities and the definition of the characteristic function
(\ref{chfn}) yield (cf. Lemmas 2.2 and 3.6 in \cite{BES1})
\begin{equation}\label{LTheta-1}L_T^* L_T =
I_{H^2_{\cld_{T^*}}(\mathbb{D})} - M_{\theta_T}
M_{\theta_T}^*,\end{equation}where $M_{\theta_T}$ is the
multiplication operator defined by $M_{\theta_T} f = \theta_T(w)
f(w)$ for all $f \in H^2_{\cld_T}(\mathbb{D})$ and $w \in
\mathbb{D}$. See \cite{BES1} for more details, where this is
carried out for a tuple of operators satisfying a {{\it ball type
condition}}.

Now we can focus on $n$ tuples of commuting operators. From this
point on, we shall assume that $n$ is an integer and $n \geq 2$. We
shall denote by $\mathbb{N}^n$ the set of all multi-indices $\bm{k}
:= (k_1, \ldots, k_n)$ where $k_i \in \mathbb{N}$ for $i = 1,
\ldots, n$. For a multi-index $\bm{k} \in \mathbb{N}^n$ we denote
${z}^{\bm{k}} = z_1^{k_1} \cdots z_n^{k_n}$ and ${T}^{\bm{k}} =
T_1^{k_1}\cdots T_n^{k_n}$ where $\bm{z} := (z_1, \ldots, z_n) \in
\mathbb{C}^n$ and ${T} = (T_1, \ldots, T_n)$ a commuting tuple (that
is, $T_i T_j = T_j T_i$ for $i, j = 1, \ldots, n$) of operators on
some Hilbert space $\clh$.

Now, we introduce the notion of isometric dilation of an $n$-tuple
operators (cf. \cite{JS}). Let $T$ and $V$ be $n$-tuples of
operators on Hilbert spaces $\clh$ and $\clk$, respectively. Then
$V$ is said to be a \textit{dilation} of $T$ if there exists an
isometry $\Pi : \clh \raro \clk$ such that
\[\Pi T_i^* = V_i^* \Pi. \quad \quad \quad (1 \leq i \leq n)\]The
dilation is said to be \textit{minimal} if \[\clk =
\overline{\mbox{span}}\{V^{\bm{k}} (\Pi \clh) : \bm{k} \in
\mathbb{N}^n\}.\]Note that $V$ on $\clk$ is a dilation of $T$ on
$\clh$ if and only if \[T_i \cong P_{\clq} V_i|_{\clq}, \quad
\quad \quad (1 \leq i \leq n)\]where $\clq$ is a joint $(V_1^*,
\ldots, V_n^*)$-invariant subspace of $\clk$ (see Section 2 of
\cite{JS} for more details).

Let $T = (T_1, \ldots, T_n)$ be an $n$-tuple of doubly commuting
contractions on $\clh$. That is, $T$ is a commuting tuple and $T_i
T_j^* = T_j^* T_i$ for $i \neq j.$ Define the defect operator
$D_{T^*}$ by
\[D_{T^*} := \prod_{i=1}^n D_{T_i^*} = \big(\prod_{i=1}^n (I_{\clh} - T_i T_i^*)\big)^{\frac{1}{2}}.
\]
and the defect space $\cld_{T^*}$ by \[\cld_{T^*} : =
\overline{\mbox{ran}} D_{T^*} = \overline{\mbox{ran}}
\prod_{i=1}^n D_{T_i^*}.\]

The \textit{Hardy space} $H^2(\mathbb{D}^n)$ over the unit polydisc
$\mathbb{D}^n$ is the Hilbert space of all holomorphic functions $f$
on $\mathbb{D}^n$ such that
\[\|f\|_{H^2(\mathbb{D}^n)} := \bigg(\mathop{\mbox{sup}}_{0 \leq r < 1}
\mathop{\int}_{\mathbb{T}^n} |f(r \bm{z})|^2
d\bm{\theta}\bigg)^{\frac{1}{2}} < \infty,\]where $d\bm{\theta}$
is the normalized Lebesgue measure on the torus $\mathbb{T}^n$,
the distinguished boundary of $\mathbb{D}^n$, and $r\bm{z} : =
(rz_1, \ldots, r z_n)$ (cf. \cite{R}, \cite{FS}). Note also that
$H^2(\mathbb{D}^n)$ is a reproducing kernel Hilbert space
\cite{Ar} corresponding to the Szego kernel $\mathbb{S} :
\mathbb{D}^n \times \mathbb{D}^n \raro \mathbb{C}$, where
\[\mathbb{S}(\z, \w) = \prod_{i=1}^n (1 - z_i \bar{w}_i)^{-1}.
\quad \quad \quad (\z, \w \in \mathbb{D}^n)\] We denote the Banach
algebra of all bounded holomorphic functions on $\mathbb{D}^n$ by
$H^{\infty}(\mathbb{D}^n)$ equipped with the supremum norm.

Given a Hilbert space $\cle$ we identify $H^2(\mathbb{D}^n)
\otimes \cle$ with $H^2_{\cle}(\mathbb{D}^n)$ via the unitary map
$\bm{z}^{\bm{k}} \otimes \eta \mapsto \bm{z}^{\bm{k}} \eta$ for
all $\bm{k} \in \mathbb{N}^n$ and $\eta \in \cle$. Moreover, it is
easy to see that the corresponding multiplication operators by the
coordinate functions are intertwined by this unitary map.

\begin{Definition}
Let $T$ be an $n$-tuple ($n > 1$) of doubly commuting contractions
on a Hilbert space $\clh$. The tuple is said to be a
\textit{doubly commuting pure tuple} if $T_i \in C_{\cdot 0}$ for
all $1 \leq i \leq n$.
\end{Definition}

The tuple of shift operators $(M_{z_1}, \ldots, M_{z_n})$ on
$H^2_{\cle}(\mathbb{D}^n)$ is a natural example of a doubly
commuting pure tuple of operators.

\newsection{Isometric dilation}

In this section we will be concerned with the isometric dilation
of a doubly commuting pure tuple on a Hilbert space $\clh$ .
Suppose that $T = (T_1, \ldots, T_n)$  is a doubly commuting
tuple. Then
\begin{equation}\label{int-1}T_i D_{T_j^*} = D_{T_j^*} T_i\end{equation}for $1 \leq i, j \leq n$ and
$i \neq j$ and \begin{equation}\label{int-2}D_{T_i^*} D_{T_j^*} =
D_{T_j^*} D_{T_i^*}. \quad \quad \quad (1 \leq i < j \leq
 n)\end{equation}

\begin{Theorem}\label{dil}
Let $T$ be a doubly commuting pure tuple on $\clh$. Then the bounded
linear operator $L_T : \clh \rightarrow
H^2_{\cld_{{T}^*}}(\mathbb{D}^n)$ defined by
\[(L_T h)(\bm{z}) = D_{T^*} \prod_{i=1}^n  (I - z_i
T^*_i)^{-1}h\]is an isometry and \[L_T T_i^* = M_{z_i}^* L_T,\]for
$i = 1, \ldots, n$. Moreover,
\[L^*_T(\mathbb{S} (\cdot, \bm{w}) \eta) = \prod_{i=1}^n (I -
\bar{w_i}T_i)^{-1} D_{T^*} \eta,\]for all $\bm{w} \in \mathbb{D}^n$
and $\eta \in \cld_{{T}^*}$, and \[H^2_{\cld_{T^*}}(\mathbb{D}^n) =
\overline{\mbox{span}}\{z^{\bm{k}} (L_T \clh) : \bm{k} \in
\mathbb{N}^n\}.\]
\end{Theorem}

\NI\textsf{Proof.} First identify $H^2_{\cld_{T_i}}(\mathbb{D}^n)$
with $H^2(\mathbb{D}) \otimes \cdots \otimes (H^2(\mathbb{D})
\otimes \cld_{T_i}) \otimes \cdots \otimes H^2(\mathbb{D})$ and
$H^2_{\cld_{{T}_i^*}}(\mathbb{D}^n)$ with $H^2(\mathbb{D})
\otimes\cdots \otimes (H^2(\mathbb{D}) \otimes \cld_{{T}_i^*})
\otimes \cdots \otimes H^2(\mathbb{D}).$ Let $T = (T_1, \ldots,
T_n)$ be a doubly commuting pure tuple on $\clh$. Then (\ref{e-1})
shows that the operator $L_{T_i} : \clh \rightarrow
{H^2_{\cld_{T_i^*}}(\mathbb{D}^n)}$ defined by
\[(L_{T_i} h)(\bm{z}) = D_{T^*_i}(I-z_iT^*_i)^{-1}h,\quad \quad \quad (h \in \clh,\,\z \in \mathbb{D}^n)\]is an
isometry for $i = 1, \ldots, n$. We now calculate
\[\begin{split}
\| h\|_{\clh}^2 & = \|L_{T_1} h\|_{H^2(\mathbb{D}^n) \otimes
\cld_{T_1^*}}^2 = \|\sum_{k_1 \in \mathbb{N}} z_1^{k_1} D_{T^*_1}
T_1^{*k_1} h\|_{H^2(\mathbb{D}^n) \otimes \cld_{T_1^*}}^2\\ &  =
\sum_{k_1 \in \mathbb{N}} \|D_{T_1^*} T_1^{*k_1} h
\|_{\cld_{T_1^*}}^2 = \sum_{k_1 \in \mathbb{N}} \|L_{T_2} (D_{T_1^*}
T_1^{*k_1} h)\|_{H^2(\mathbb{D}^n) \otimes \cld_{T_2^*}}^2 \notag\\
& = \sum_{k_1 \in \mathbb{N}} \| \sum_{k_2 \in \mathbb{N}} z_2^{k_2}
D_{T_2^*} T_2^{*k_2} D_{T_1^*} T_1^{*k_1} h\|_{H^2(\mathbb{D}^n)
\otimes \cld_{T_2^*}}^2 = \sum_{k_1, k_2 \in \mathbb{N}}
\|D_{T^*_2}D_{T^*_1} T_1^{*k_1} T_2^{*k_2} h\|_{\cld_{T_2^*}}^2\\ &
= \sum_{k_1, k_2 \in \mathbb{N}} \|D_{T^*_1}D_{T^*_2} T_1^{*k_1}
T_2^{*k_2} h\|_{\overline{\mbox{ran}} (D_{T_1^*} D_{T_2^*})}^2.
\quad \quad \quad (h \in \clh)
\end{split}\] Continuing this process we obtain
\[\|h\|_{\clh}^2 = \sum_{\bm{k} \in \mathbb{N}^n} \|\prod_{i=1}^n D_{T^*_i}
{T}^{*\bm{k}} h\|_{\overline{\mbox{ran}} (D_{T_1^*} \cdots
D_{T_n^*})}^2 = \sum_{\bm{k} \in \mathbb{N}^n} \|D_{T^*}
{T}^{*\bm{k}} h\|^2_{\cld_{T^*}}.\]Hence it follows that
\[\|h\|_{\clh}^2 = \| \sum_{\bm{k} \in
\mathbb{N}^n} z^{\bm{k}} D_{T^*} T^{*\bm{k}}
h\|^2_{H^2_{\cld_{T^*}}(\mathbb{D}^n)} = \|L_T
h\|^2_{H^2_{\cld_{T^*}}(\mathbb{D}^n)}. \quad \quad \quad (h \in
\clh)\]This implies that $L_T$ is an isometry. Moreover \[L_{T}
T_i^* h = D_{T^*} \sum_{\bm{k} \in \mathbb{N}^n} z^{\bm{k}}
T^{*(\bm{k} + e_i)} h = M^*_{z_i} D_{T^*} \sum_{\bm{k} \in
\mathbb{N}^n} z^{\bm{k}} T^{* \bm{k} } h = M_{z_i}^* L_{T} h, \quad
\quad (h \in \clh, 1 \leq i \leq n)\]and consequently \[L_{T} T_i^*
= M_{z_i}^* L_{T}. \quad \quad \quad (1 \leq i \leq n)\] Also for
all $h \in \clh$, $\eta \in \cld_{T^*}$ and $\w \in \mathbb{D}^n$,
it follows that
\[\begin{split}\langle L^*_{{T}}(\mathbb{S}(\cdot, \bm{w}) \eta), h\rangle_{\clh} & = \langle \mathbb{S}(\cdot, \bm{w}) \eta,
L_{\bm{T}} h \rangle_{H^2_{\cld_{T^*}}(\mathbb{D}^n)} \\ & = \langle
\sum_{\bm{k} \in \mathbb{N}^n} {z}^{\bm{k}} \bar{{w}}^{\bm{k}} \eta,
\sum_{\bm{l} \in
\mathbb{N}^n} {z}^{\bm{l}} D_{T^*} {T}^{*\bm{l}} h\rangle_{H^2_{\cld_{T^*}}(\mathbb{D}^n)} \\
 & = \sum_{\bm{k} \in
\mathbb{N}^n} \langle \bar{{w}}^{\bm{k}} \eta, D_{T^*} {T}^{*\bm{k}}
h \rangle_{\clh},\end{split}\]and so
\[\langle L^*_{{T}}(\mathbb{S} (\cdot, \bm{w}) \eta), h\rangle_{\clh}  = \langle \prod_{i=1}^n (I - \bar{w}T_i)^{-1}
D_{T^*} \eta, h\rangle_{\clh}.
\]

\NI We complete the proof by showing that the dilation $(M_{z_1},
\ldots, M_{z_n})$ on $H^2_{\cld_{T^*}}(\mathbb{D}^n)$ is minimal,
that is,\[H^2_{\cld_{T^*}}(\mathbb{D}^n) =
\overline{\mbox{span}}\{z^{\bm{k}} (L_T \clh) : \bm{k} \in
\mathbb{N}^n\}.\]But since $\overline{\mbox{span}}\{z^{\bm{k}} (L_T
\clh) : \bm{k} \in \mathbb{N}^n\}$ is a joint $(M_{z_1}, \ldots,
M_{z_n})$-reducing closed subspace of
$H^2_{\cld_{T^*}}(\mathbb{D}^n)$, it follows from Proposition 2.2 in
\cite{SSW} that
\[\overline{\mbox{span}}\{z^{\bm{k}} (L_T \clh) : \bm{k} \in
\mathbb{N}^n\} = H^2_{\cle}(\mathbb{D}^n),\]for some $\cle \subseteq
\cld_{T^*}$. We claim that $\cle = \cld_{T^*}$. To see that, first
we note that for $(M_{z_1}, \ldots, M_{z_n})$ on
$H^2_{\cld_{T^*}}(\mathbb{D}^n)$ we have (cf. \cite{SSW})
\[\mathop{\sum}_{0 \leq i_1 < \ldots < i_l \leq n} (-1)^l M_{z_{i_1}}
\cdots M_{z_{i_l}} M^*_{{z}_{i_1}} \cdots M^*_{{z}_{i_l}} =
P_{\cld_{T^*}},\] where $P_{\cld_{T^*}}$ is the projection to the
space of constant functions. We then have
\[\big(\mathop{\sum}_{0 \leq i_1 < \ldots < i_l \leq n} (-1)^l
M_{z_{i_1}} \cdots M_{z_{i_l}} M^*_{{z}_{i_1}} \cdots
M^*_{{z}_{i_l}}\big) (L_T h) = P_{\cld_{T^*}}(L_T h) = (L_T h)(0).
\quad \quad (h \in \clh)\]On the other hand,
\[(L_T h)(0) = (D_{T^*} \prod_{i=1}^n (I - z_i T_i^*)^{-1} h)(0) =
D_{T^*}h.\]It now follows that $\cle = \cld_{T^*}$ and the proof is
complete. \qed

The following corollary is a rephrasing of the definition of
isometric dilation and Theorem \ref{dil}.

\begin{Corollary}\label{Q-dc}
Let $T$ be a doubly commuting pure tuple on $\clh$. Then
$(M_{z_1},\ldots, M_{z_n})$ on $H^2_{\cld_{T^*}}(\mathbb{D}^n)$ is
the minimal isometric dilation of $T$, that is, there exists a joint
$(M_{z_1}^*, \ldots, M_{z_n}^*)$-invariant subspace $\clq$ of
$H^2_{\cld_{{T}^*}}(\mathbb{D}^n)$ such that
\[T_i \cong P_{\clq} M_{z_i}|_{\clq},\]for all $1 \leq i \leq n$, and
\[H^2_{\cld_{T^*}}(\mathbb{D}^n) =
\overline{\mbox{span}}\{z^{\bm{k}} \clq : \bm{k} \in
\mathbb{N}^n\}.\]

\end{Corollary}

The proofs of the dilation theorem obtained in this way are quite
different from any earlier proofs (cf. \cite{MV}, \cite{CV},
\cite{sub3}, \cite{JS}).

\bigskip
\noindent{\bf Remark:} An anonymous referee of an earlier version
of this paper pointed out that most of Theorem 3.1 can be obtained
using results from \cite{AEM}. But, our proofs are essentially
arguments based on the case of a single contraction (unlike that
of \cite{AEM}), because the deflect operator splits into a product
of individual defect operators. Hence, our techniques demonstrate
the importance of the dilation theory of a single contraction in
the dilation theory of a tuple of doubly commuting contractions.

\newsection{Canonical model}

In this section, we study the analytic structure of the backward
shift invariant subspace $\clq$ in Corollary \ref{Q-dc}. We begin
with a few definitions.

Let $T = (T_1, \ldots, T_n)$ be an $n$-tuple of commuting
contractions on $\clh$. Define a one variable multiplier
$\Theta_{T_i} \in H^\infty_{\clb(\cld_{T_i},
\cld_{T^*_i})}(\mathbb{D}^n)$ by
\[\Theta_{T_i}(\z) = \theta_{T_i}(z_i), \quad \quad \quad (\z \in
\mathbb{D}^n)\]where $\theta_{T_i}$ is the characteristic function
of the contraction $T_i$ and $i = 1, \ldots, n$ (see the definition
in (\ref{chfn})). Therefore, $M_{\Theta_{T_i}} :
H^2_{\cld_{T_i}}(\mathbb{D}^n) \raro
H^2_{\cld_{T_i^*}}(\mathbb{D}^n)$ is a bounded linear operator
defined by
\begin{equation}\label{Theta-1}(M_{\Theta_{T_i}} f)(\z) = (\Theta_{T_i} f)(\bm{z}) =
\theta_{T_i}(z_i) f(\bm{z}), \quad \quad \quad (\bm{z} \in
\mathbb{D}^n, f \in H^2_{\cld_{T_i}}(\mathbb{D}^n))\end{equation}
for $i = 1, \ldots, n$. It is easy to see that \[M_{\Theta_{T_i}}
M_{z_j} = M_{z_j} M_{\Theta_{T_i}},\]for all $i, j = 1, \ldots, n$,
and
\[M_{\Theta_{T_i}} M^*_{z_j} = M_{z_j}^* M_{\Theta_{T_i}},\]for all $i,
j = 1, \ldots, n$, and $i \neq j$. We have,
\[M_{\Theta_{T_i}} = I_{H^2(\mathbb{D})} \otimes \cdots \otimes M_{\theta_{T_i}}
\otimes \cdots \otimes I_{ H^2(\mathbb{D})},\]for $i = 1, \ldots,
n$. We have also by virtue of (\ref{e-3})
\begin{equation}\label{4.1}\mathbb{S}(z_i, w_i)(I_{\cld_{T_i^*}} -
\Theta_{T_i}(\z) \Theta_{T_i}(\w)^*) = D_{T_i^*} (I - z_i
T_i^*)^{-1} (I - \bar{w}_i T_i)^{-1} D_{T_i^*},\end{equation}for
$i = 1, \ldots, n$.  Here we record the following simple
observation.

\begin{Lemma}
Let ${T} = (T_1, \ldots, T_n)$ be an $n$-tuple of bounded linear
operators on a Hilbert space $\clh$ and $T_i$ be of class
$C_{\cdot 0}$ for all $1 = 1, \ldots, n$. Then $\Theta_{T_i} \in
H^{\infty}_{\clb(\cld_{T_i}, \cld_{{T}^*})}(\mathbb{D}^n)$ is a
one variable inner function for $i = 1, \ldots, n$.
\end{Lemma}

Now suppose that $T$ is a doubly commuting tuple. The equalities
(\ref{int-1}) and (\ref{int-2}) imply that
\[[D_{T_i^*} (I - z_i T_i^*)^{-1} (I - \bar{w}_i T_i)^{-1}
D_{T_i^*}] \big( \prod_{j=1}^n D_{T^*_j}\big) = \big( \prod_{j=1}^n
D_{T^*_j}\big) [(I - z_i T_i^*)^{-1} (I - \bar{w}_i T_i)^{-1}
D_{T_i^*}^2] ,\]and hence
\begin{equation}\label{4a}[D_{T_i^*} (I - z_i T_i^*)^{-1} (I -
\bar{w}_i T_i)^{-1} D_{T_i^*}]\cld_{T^*} \subseteq \cld_{T^*},
\end{equation} for $i = 1, \ldots, n$. This observation, together
with (\ref{4.1}) imply that
\begin{equation}\label{4b}(\Theta_{T_i}(\z) \Theta_{T_i}(\w)^*)
\cld_{T^*} \subseteq \cld_{T^*}.\end{equation}In particular,
\begin{equation}\label{T-inclusion}(M_{\Theta_{T_i}} M_{\Theta_{T_i}}^*)
H^2_{\cld_{T^*}}(\mathbb{D}^n) \subseteq
H^2_{\cld_{T^*}}(\mathbb{D}^n). \quad \quad \quad (1 \leq i \leq
n)\end{equation} Moreover, it follows from (\ref{4.1}), (\ref{4a})
and (\ref{4b}) that
\begin{equation}\label{E2}\prod_{i=1}^n[D_{T_i^*} (I - z_i
T_i^*)^{-1} (I - \bar{w}_i T_i)^{-1} D_{T_i^*}]|_{\cld_{T^*}} =
\mathbb{S}(\z, \w) \prod_{i=1}^n (I_{\cld_{T_i^*}} -
\Theta_{T_i}(\z) \Theta_{T_i}(\w)^*)|_{\cld_{T^*}}.\end{equation}

The following result relates the characteristic functions of the
coordinate operators and the isometric dilation of a doubly
commuting pure tuple $T$.

\begin{Proposition}\label{L-theta-n} Let ${T}$ be a doubly commuting
pure tuple of operators on $\clh$. Then
\[L_{T} L_{{T}}^* = \prod_{i=1}^n (I_{H^2_{\cld_{T^*}}(\mathbb{D}^n)}
- M_{\Theta_{T_i}} M_{\Theta_{T_i}}^*|_{H^2_{\cld_{T^*}(\mathbb{D}^n)}}).\]
\end{Proposition}
\NI\textsf{Proof.} Let $\z, \w \in \mathbb{D}^n$ and $\eta, \zeta
\in \cld_{{T}^*}$ so that
\[\begin{split} \langle L_{{T}} L_{{T}}^*(\mathbb{S}(\cdot, \w)
\eta),\, \mathbb{S} (\cdot, \z) \zeta
\rangle_{H^2_{\cld_{T^*}(\mathbb{D}^n)}} & = \langle \prod_{i=1}^n
(I - \bar{w}_i T_i)^{-1} D_{T^*} \eta, \prod_{j=1}^n (I - \bar{z}_j
T_j)^{-1} D_{T^*} \zeta \rangle_{\clh} \notag \\ & = \langle
\prod_{i=1}^n D_{T^*} (1 - z_i T^*_i)^{-1} (I - \bar{w}_i T_i)^{-1}
D_{T^*} \eta, \zeta\rangle_{\clh}.\end{split}\] By virtue of
(\ref{E2}), it follows that
\[\begin{split} \langle L_{{T}} L_{{T}}^*(\mathbb{S}(\cdot, \w)
\eta), \mathbb{S}(\cdot, \z) \eta
\rangle_{H^2_{\cld_{T^*}(\mathbb{D}^n)}} & = \mathbb{S}(\z, \w)
\langle \prod_{i=1}^n(I - \Theta_{T_i}(\z) \Theta_{T_i}(\w)^*) \eta,
\zeta \rangle \notag
\\& = \langle \prod_{i=1}^n (I_{H^2_{\cld_{T^*}}(\mathbb{D}^n)} -
M_{\Theta_{T_i}} M_{\Theta_{T_i}}^*)(\mathbb{S} (\cdot, \w) \eta),
\mathbb{S}(\cdot, \z) \eta \rangle,
\end{split}\]
which completes the proof of the proposition. \qed

The following well known result (cf. \cite{JS2}), concerning the
range of the sum of a finite family of commuting orthogonal
projections, will play a key role in the model theory for doubly
commuting pure tuples.

\begin{Lemma}\label{P-F} Let $\{P_i\}_{i=1}^n$ be a collection of commuting orthogonal
projections on a Hilbert space $\clh$. Then $\cll :=
\mathop{\sum}_{i=1}^n \mbox{ran} P_i$ is closed and the orthogonal
projection of $\clh$ onto $\cll$ is given by
\[P_{\cll}
= I_{\clh} - \mathop{\prod}_{i=1}^n (I_{\clh} - P_i).\]
\end{Lemma}
\NI\textsf{Proof.} We set $X_i = P_i (I_{\clh} - P_{i+1}) \cdots
(I_{\clh} - P_{n-1}) (I_{\clh} - P_n)$ for all $i = 1, \ldots, n-1$,
and $X_n = P_n$. Since \[\mathop{\sum}_{i=1}^n X_i = I_{\clh} -
\mathop{\Pi}_{i=1}^n (I_{\clh} - P_i),\]and $\{X_i\}_{i=1}^n$ is a
family of orthogonal projections with orthogonal ranges, we have
\[\cll = \mbox{ran} X_1 \oplus \cdots \oplus
\mbox{ran} X_n.\] This completes the proof of the lemma. \qed

We now have the following key corollary to the main result of this
section.

\begin{Corollary}\label{key}
Let $T$ be a doubly commuting pure tuple on $\clh$. Then \[\cls_T :
= \sum_{i=1}^n \big(H^2_{\cld_{T^*}}(\mathbb{D}^n) \bigcap
\Theta_{T_i} H^2_{\cld_{T_i}}(\mathbb{D}^n)\big)\]is a closed
subspace of $H^2_{\cld_{T^*}}(\mathbb{D}^n)$ and
\[I_{H^2_{\cld_{T^*}}(\mathbb{D}^n)} - P_{\cls_T} =
\prod_{i=1}^n (I_{H^2_{\cld_{T^*}}(\mathbb{D}^n)} - M_{\Theta_{T_i}}
M_{\Theta_{T_i}}^*)|_{H^2_{\cld_{T^*}}(\mathbb{D}^n)}.\]
\end{Corollary}
\NI\textsf{Proof.} It follows from the definition of
$M_{\Theta_{T_i}}$ and the fact that $T_i$ is pure, that
$M_{\Theta_{T_i}}$ is an isometry and hence $M_{\Theta_{T_i}}
M_{\Theta_{T_i}}^*$ is an orthogonal projection for $i = 1,
\ldots, n$. Also by (\ref{T-inclusion}), we have
\[P_{H^2_{\cld_{T^*}}(\mathbb{D}^n)}(M_{\Theta_{T_i}}
M_{\Theta_{T_i}}^*)P_{H^2_{\cld_{T^*}}(\mathbb{D}^n)} =
(M_{\Theta_{T_i}}
M_{\Theta_{T_i}}^*)P_{H^2_{\cld_{T^*}}(\mathbb{D}^n)}.\] Let $P_i =
(M_{\Theta_{T_i}}
M_{\Theta_{T_i}}^*)|_{H^2_{\cld_{T^*}}(\mathbb{D}^n)} \in
\clb(H^2_{\cld_{T^*}}(\mathbb{D}^n))$. Then $P_i$, for each $i = 1,
\ldots, n$, is an orthogonal projection and
\begin{equation}\label{theta-proj}\mbox{ran} P_i = \mbox{ran} M_{\Theta_{T_i}} \bigcap
H^2_{\cld_{T^*}}(\mathbb{D}^n) = {\Theta_{T_i}}
H^2_{\cld_{T_i}}(\mathbb{D}^n)\bigcap
H^2_{\cld_{T^*}}(\mathbb{D}^n).\end{equation}Further,
\[P_i P_j = P_j P_i. \quad \quad (1 \leq i < j \leq
n)\]By Lemma \ref{P-F} and (\ref{theta-proj}), we have \[\cls_T =
\sum_{i=1}^n \mbox{ran} P_i = \sum_{i=1}^n
\big(H^2_{\cld_{T^*}}(\mathbb{D}^n) \bigcap \Theta_{T_i}
H^2_{\cld_{T_i}}(\mathbb{D}^n)\big),\]is a closed subspace of
$H^2_{\cld_{T^*}}(\mathbb{D}^n)$. Again by Lemma \ref{P-F}, we have
\[P_{\cls_T} = I_{H^2_{\cld_{T^*}}(\mathbb{D}^n)} - \prod_{i=1}^n
(I_{H^2_{\cld_{T^*}}(\mathbb{D}^n)} - P_i) =
I_{H^2_{\cld_{T^*}}(\mathbb{D}^n)} - \prod_{i=1}^n
(I_{H^2_{\cld_{T^*}}(\mathbb{D}^n)} - M_{\Theta_{T_i}}
M_{\Theta_{T_i}}^*)|_{H^2_{\cld_{T^*}}(\mathbb{D}^n)}.\]This
completes the proof. \qed

\begin{Theorem}\label{M-1}
Let $T$ be a doubly commuting pure tuple on $\clh$. Then for all $i
= 1, \ldots, n$, \[T_i \cong P_{\clq_{{T}}}
M_{z_i}|_{\clq_{{T}}},\]where
\[\clq_{{T}} = \cls_T^\perp \cong H^2_{\cld_{{T}^*}}(\mathbb{D}^n) /\cls_T,\]
is a joint $(M_{z_1}^*, \ldots, M_{z_n}^*)$-invariant subspace of
$H^2_{\cld_{{T}^*}}(\mathbb{D}^n)$ corresponding to the joint
$(M_{z_1}, \ldots, M_{z_n})$-invariant subspace
\[\cls_T =  \sum_{i=1}^n \big(H^2_{\cld_{T^*}}(\mathbb{D}^n) \bigcap
\Theta_{T_i} H^2_{\cld_{T_i}}(\mathbb{D}^n)\big).\]
\end{Theorem}
\NI\textsf{Proof.} Let ${T}$ be a doubly commuting pure tuple on
$\clh$. By Proposition \ref{L-theta-n}, we have
\[L_{{T}} L_{{T}}^* = \prod_{i=1}^n
(I_{H^2_{\cld_{\bm{T}^*}}(\mathbb{D}^n)} - M_{\Theta_{T_i}}
M_{\Theta_{T_i}}^*)|_{H^2_{\cld_{{T}^*}}(\mathbb{D}^n)}.\]This along
with Corollary \ref{key} yields
\[ \begin{split} L_{T} L_{{T}}^* & = I_{H^2_{\cld_{T^*}}(\mathbb{D}^n)} - [I_{H^2_{\cld_{T^*}}(\mathbb{D}^n)} - \prod_{i=1}^n
(I_{H^2_{\cld_{\bm{T}^*}}(\mathbb{D}^n)} - M_{\Theta_{T_i}}
M_{\Theta_{T_i}}^*)]|_{H^2_{\cld_{{T}^*}}(\mathbb{D}^n)}\notag\\
& = I_{H^2_{\cld_{{T}^*}}(\mathbb{D}^n)} -
P_{\cls_T}.\end{split}\]Consequently, \[\mbox{ran} L_{{T}} \cong
\cls_T^{\perp} \cong H^2_{\cld_{{T}^*}}(\mathbb{D}^n) / \cls_T,
\]and \[T_i \cong P_{\clq_{{T}}} M_{z_i}|_{\clq_{{T}}},\]for
$i = 1, \ldots, n$. This completes the proof. \qed

\section{One variable inner functions}

The purpose of this section is to obtain a concrete realization of
the joint $(M_{z_1}, \ldots, M_{z_n})$-invariant subspace $\cls_T$,
in Theorem \ref{M-1}, in terms of one variable inner functions on
the polydisc.

Let $T$ be a doubly commuting pure tuple of operators on $\clh$. By
Theorem \ref{M-1}, we get \[\clh \cong \cls_T^\perp,\quad
\mbox{and}\quad T_i \cong P_{\cls_T^\perp}
M_{z_i}|_{\cls_T^\perp},\]where \[\cls_T = \sum_{i=1}^n
\cls_{T_i},\] is a joint $(M_{z_1}, \ldots, M_{z_n})$-invariant
subspace of $H^2_{\cld_{T^*}}(\mathbb{D}^n)$ and
\[\cls_{T_i} := H^2_{\cld_{T^*}}(\mathbb{D}^n) \bigcap \Theta_{T_i}
H^2_{\cld_{T_i}}(\mathbb{D}^n).\quad \quad \quad (1 \leq i \leq n)\]
Recall that $H^2_{\cld_{T^*}}(\mathbb{D}^n)$ and $\Theta_{T_i}
H^2_{\cld_{T_i}}(\mathbb{D}^n)$ can be identified with
$H^2(\mathbb{D}) \otimes \cdots \otimes H^2_{\cld_{T^*}}(\mathbb{D})
\otimes \cdots \otimes H^2(\mathbb{D})$ and $H^2(\mathbb{D}) \otimes
\cdots \otimes \big(\theta_{T_i} H^2_{\cld_{T_i}}(\mathbb{D})\big)
\otimes \cdots \otimes H^2(\mathbb{D})$, respectively. Also
\[\cls_{T_i} \cong H^2(\mathbb{D}) \otimes \cdots \otimes
\tilde{\cls}_{T_i} \otimes \cdots \otimes H^2(\mathbb{D}),\]for
some $M_z$-invariant subspace $\tilde{\cls}_{T_i}$ of
$H^2_{\cld_{T_i^*}}(\mathbb{D})$.

\NI Let $1 \leq i \leq n$ and assume that $\cls_{T_i} \neq \{0\}$.
Then by the Beurling-Lax-Halmos theorem, on shift invariant
subspaces of vector-valued Hardy spaces (\cite{SNF}), there exists
a Hilbert space $\cle_{T_i}$ and an inner multiplier $\phi_{T_i}
\in H^\infty_{\clb(\cle_{T_i}, \cld_{T^*})}(\mathbb{D})$, unique
up to unitary equivalence, such that
\[\tilde{\cls}_{T_i} = \phi_{T_i}
H^2_{\cle_{T_i}}(\mathbb{D}).\]Thus\[\cls_{T_i} \cong
H^2(\mathbb{D}) \otimes \cdots \otimes \big(\phi_{T_i}
H^2_{\cle_{T_i}}(\mathbb{D})\big) \otimes \cdots \otimes
H^2(\mathbb{D}).\] Let  \[(\Phi_{T_i} f)(\z) = \phi_{T_i}(z_i)
f(\z). \quad \quad \quad (\z \in \mathbb{D}^n, f \in
H^2_{\cle_{T_i}}(\mathbb{D}^n))\]Certainly $\Phi_{T_i} \in
H^\infty_{\clb(\cle_{T_i}, \cld_{T^*})}(\mathbb{D}^n)$ is a one
variable inner function. Moreover, $H^2(\mathbb{D}) \otimes \cdots
\otimes \big(\phi_{T_i} H^2_{\cle_{T_i}}(\mathbb{D})\big) \otimes
\cdots \otimes H^2(\mathbb{D})$ can be identified to $\Phi_{T_i}
H^2_{\cle_{T_i}}(\mathbb{D}^n)$, via the same identification map,
and
\[\tilde{\cls}_{T_i} = \Phi_{T_i}
H^2_{\cle_{T_i}}(\mathbb{D}^n).\]Consequently, \[\cls_{T} =
\sum_{i=1}^n \Phi_{T_i} H^2_{\cle_{T_i}}(\mathbb{D}^n),\]where
each $\Phi_{T_i} \in H^\infty_{\clb(\cle_{T_i},
\cld_{T^*})}(\mathbb{D}^n)$ is either a one variable inner
function in $z_i$, or the zero function and $i= 1, \ldots, n$.

This along with Theorem  \ref{M-1} proves the following result.

\begin{Theorem}\label{m-2}
Let $T$ be a doubly commuting pure tuple on $\clh$. Then there
exists a joint $(M_{z_1}^*, \ldots, M_{z_n}^*)$-invariant subspace
$\clq_T$ of $H^2_{\cld_{T^*}}(\mathbb{D}^n)$ such that \[\clh
\cong \clq_T, \quad \mbox{and} \quad T_i \cong P_{\clq_T}
M_{z_i}|_{\clq_T},\]for $i = 1, \ldots, n$. Moreover, there exists
Hilbert spaces $\{\cle_{T_i}\}_{i=1}^n$ and $\Phi_{T_i} \in
H^\infty_{\clb(\cle_{T_i}, \cld_{T^*})}(\mathbb{D}^n)$, unique up
to unitary equivalence, such that each $\Phi_{T_i}$ ($1 \leq i
\leq n$) is either a one variable inner function in $z_i$, or the
zero function and
\[\cls_{T} := \sum_{i=1}^n \Phi_{T_i}
H^2_{\cle_{T_i}}(\mathbb{D}^n) \]is closed in
$H^2_{\cld_{T^*}}(\mathbb{D}^n)$, and \[\clq_T = \cls_T^\perp.\]
\end{Theorem}
In particular, Theorem \ref{m-2} says that the class of all doubly
commuting pure tuples on separable Hilbert spaces is equal, up to
unitary equivalence, to the class of all doubly commuting
$(M_{z_1}^*, \ldots, M_{z_n}^*)$-invariant subspaces of
vector-valued Hardy spaces over the polydisc.

As a special case of Theorem \ref{m-2} we obtain the following
corollary.

\begin{Corollary}
Let $\clq$ be a joint $(M_{z_1}^*, \ldots, M_{z_n}^*)$-invariant
closed proper subspace of $H^2(\mathbb{D}^n)$ and let $C_{z_i}: =
P_{\clq} M_{z_i}|_{\clq}$ for $i= 1, \ldots, n$. Then $(C_{z_1},
\ldots, C_{z_n})$ is doubly commuting if and only if there exists
$\{\theta_i\}_{i=1}^n \subseteq H^\infty(\mathbb{D})$ such that each
$\theta_i$ is either inner or the zero function for $i = 1, \ldots,
n$ and \[\clq = \big(\sum_{i=1}^n \Theta_i
H^2(\mathbb{D}^n)\big)^\perp,\]where $\Theta_i(\z) = \theta_i(z_i)$
for all $\z \in \mathbb{D}^n$ and $i = 1, \ldots, n$.
\end{Corollary}

\NI\textsf{Proof.} If $T := (C_{z_1}, \ldots, C_{z_n})$, then
\[\begin{split}D_{T^*}^2 & = \prod_{i=1}^n (I_{\clq} - C_{z_i}
C_{z_i}^*) = P_{\clq} \big(\prod_{i=1}^n (I_{H^2(\mathbb{D}^n)} -
M_{z_i} M_{z_i}^*)\big)|_{\clq} = P_{\clq}
P_{\mathbb{C}}|_{\clq}.\end{split}\]Thus the rank of $D_{T^*}$ is
one. Now the result follows from Theorem \ref{m-2}. \qed

This result was proved by the third author in \cite{JS2}. See also
the work by Izuchi, Nakazi and Seto \cite{INS} for the base case $n
= 2$.

\end{document}